\documentclass{article}
\usepackage{amsmath,amssymb,amsthm}
\usepackage{color}
\usepackage{hyperref}

\newtheorem{theorem}{Theorem}
\newtheorem{proposition}{Proposition}
\newtheorem{lemma}[theorem]{Lemma}
\theoremstyle{definition}
\newtheorem{example}[theorem]{Example}

\usepackage{todonotes}

\newcommand{\PAF}{\operatorname{PAF}}
\newcommand{\DFT}{\operatorname{DFT}}
\newcommand{\PSD}{\operatorname{PSD}}
\newcommand{\Fix}{\operatorname{Fix}}
\newcommand{\Aut}{\operatorname{Aut}}
\newcommand{\lcm}{\operatorname{lcm}}

\newcommand{\Q}{\mathbb{Q}}

\begin{document}

\title{Quaternary Legendre pairs II}

\author{Ilias S.~Kotsireas$^1$
  \and Christoph Koutschan$^2$
  \and Arne Winterhof $^2$}
\date{{\normalsize
  $^1$CARGO Lab, Wilfrid Laurier University, Waterloo, ON N2L 3C5, Canada\\
  $^2$RICAM, Austrian Academy of Sciences, Altenberger Stra\ss e  69, 4040~Linz, Austria\\
  e-mail: ikotsire@wlu.ca, christoph.koutschan@ricam.oeaw.ac.at,
  arne.winterhof@oeaw.ac.at}\\[2ex]
  \today}

\maketitle

\begin{abstract}
Quaternary Legendre pairs are pertinent to the construction of quaternary Hadamard matrices and have many applications, for example in coding theory and communications.

In contrast to binary Legendre pairs, quaternary ones can exist for even length $\ell$ as well.
It is conjectured that there is a quaternary Legendre pair for any even $\ell$. The smallest open case until now had been $\ell=28$, and $\ell=38$ was the only length $\ell$ with $28\le \ell\le 60$ resolved before. Here we provide constructions for $\ell=28,30,32$, and $34$. 
In parallel and independently, Jedwab and Pender found a construction of quaternary Legendre pairs of length $\ell=(q-1)/2$ for any prime power $q\equiv 1\bmod 4$, which in particular covers $\ell=30$, $36$, and $40$, so that now $\ell=42$ is the smallest unresolved case.

The main new idea of this paper is a way to separate the search for the subsequences along even and odd indices which substantially reduces the complexity of the search algorithm. 

In addition, we use Galois theory for cyclotomic fields to derive conditions which improve the PSD test.
\end{abstract}


\section{Introduction}

Two sequences 
$$A=[a_0,a_1,\ldots,a_{\ell-1}],\quad B=[b_0,b_1,\ldots,b_{\ell-1}]\in \mathbb{C}^\ell
$$ 
of period $\ell$ form a {\em Legendre pair} $(A,B)$ if
\begin{equation}\label{Legpairdef}\PAF(A,s)+\PAF(B,s)=-2,\quad s=1,2,\ldots,\left\lfloor\frac{\ell}{2}\right\rfloor,
\end{equation}
where 
$$\PAF(A,s)=\sum_{j=0}^{\ell-1} a_j\overline{a_{j+s}},\quad s=0,1,\ldots,\ell-1,$$
is the {\em periodic autocorrelation function} $\PAF(A,s)$ of $A$ at lag $s$.
(Note that $(\ref{Legpairdef})$ holds also true for $s=\left\lfloor\frac{\ell}{2}\right\rfloor+1,\ldots,\ell-1$
since $\PAF(A,\ell-s)=\overline{\PAF(A,s)}$ for $s=1,2,\ldots,\ell-1$.)

{\em Binary} Legendre pairs $(A,B)$ where $A,B \in \{-1,+1\}^{\ell}$ can be traced back to~\cite{fgs01} and \cite{s69}, see~\cite{kkbatr23,kw24} for the state of the art. In particular, the length of a binary Legendre pair must be odd, there are several known infinite classes and it is conjectured that a binary Legendre pair of length~$\ell$ exists for any odd~$\ell$. The smallest unsettled case is $\ell=115$.

Very recently~\cite{kw24}, we started studying {\em quaternary} Legendre pairs~$(A,B)$, where
$$A,B\in \{-1,+1,-i,+i\}^\ell.$$
Quaternary Legendre pairs (qLPs) are pertinent to the construction of quaternary Hadamard matrices, see \cite[Theorem~2.2]{kw24}.
In contrast to the binary case, quaternary Legendre sequences exist also for even $\ell$. 
Since a binary sequence can be considered quaternary, we may restrict ourselves to even~$\ell$.
In our previous paper~\cite{kw24} we constructed quaternary Legendre sequences for all even $\ell\le 26$ as well as for $\ell=38,62,74,82$. Note that the latter three lengths are of the form
$\ell=2p$ with a prime such that $2p-1$ is a sum of two squares, where a special \emph{compressed} pair can be used as a starting point (for the concept of compression, see~\cite{kw24}). However, this idea cannot be used for any $\ell$ of a different form including $\ell=28,30,32$ and $34$.

In parallel and independently, Jedwab and Pender \cite{jp24} constructed quaternary Legendre pairs of length $\ell=(q-1)/2$ for every prime power $q\equiv 1\bmod 4$ as well as of length $\ell=2p$ for every prime $p>2$ for which $2p-1$ is a prime power. This covers $\ell=30$ but not $\ell=28$, $32$ and $34$.

Binary and quaternary Hadamard matrices and thus binary and quaternary Legendre pairs have many applications, for example in coding theory and communication, see for example~\cite{h07}.

In this paper we continue the search for quaternary Legendre pairs of even length~$\ell$. 
A new idea for any even~$\ell$ is to separate the search for $[a_0,a_2,\ldots,a_{\ell-2}]$ and $[a_1,a_3,\ldots,a_{\ell-1}]$
or $[b_0,b_2,\ldots,b_{\ell-2}]$ and $[b_1,b_3,\ldots,b_{\ell-1}]$, respectively, which dramatically reduces the complexity of the search algorithm, see Section~\ref{general}.
In the case that $\ell$ is divisible by~$4$ we provide some additional speed-ups.


In addition, we also provide some further tests in Section~\ref{prod} based on Galois theory of cyclotomic fields and a well-known condition  
for integers which are the square of the absolute value of an algebraic integer in a cyclotomic field.

We start with some preliminary results in Section~\ref{prel}. 
Our new constructions for $\ell=28,30,32$ and $34$ are given in Section~\ref{constructions}.
(Note that our construction of length $30$ is different from the construction of \cite{jp24}.)
Now the smallest open case is $\ell=42$ and there are only $14$ unresolved cases $\le 100$:
$$\ell=42,46,52,58,64,66,70,72,76,80,88,92,94,100.$$

\section{Preliminaries}\label{prel}
\subsection{Quaternary Legendre pairs}
Let $A=[a_0,a_1,\ldots,a_{\ell-1}]$, $B=[b_0,b_1,\ldots,b_{\ell-1}]$
be a quaternary Legendre pair of even length $\ell$.

\subsubsection*{Balance}

Put
$$\alpha=\sum_{j=0}^{\ell-1}a_j \quad \mbox{and} \quad  \beta=\sum_{j=0}^{\ell-1}b_j.$$
By \cite[Lemma~1]{kw24} we may assume w.l.o.g.\ 
\begin{equation}\label{balance}\alpha=0\quad \mbox{and}\quad \beta=1+i.
\end{equation}

\subsubsection*{DFT and the PSD-test}

For an integer $n\ge 1$ we denote by $\xi_n$ the primitive $n$th root of unity 
$$\xi_n=e^{2\pi i/n}.$$
Note that we may fix any other primitive $n$th root of unity instead.

The {\em discrete Fourier transform} $\DFT(A,s)$ of $A$ at $s$
is
$$\DFT(A,s)=\sum_{j=0}^{\ell-1}a_j \xi_\ell^{js},\quad s=0,1,\ldots,\ell-1,$$
and the {\em power spectral density} $\PSD(A,s)$ of $A$ at $s$ is
$$\PSD(A,s)=|\DFT(A,s)|^2=\sum_{j=0}^{\ell-1}\PAF(A,j)\xi_\ell^{js}.$$
A quaternary Legendre pair $(A,B)$ satisfies
\begin{equation}\label{sumPSD}
  \PSD(A,s)+\PSD(B,s)=2\ell +2,\quad s=1,2,\ldots,\ell-1.
\end{equation}
For efficiency, we search for possible candidate sequences $A$ and $B$ separately.
An important tool for this separation is the following {\em PSD-test}. 
The nonnegativity of the PSD and 
\eqref{sumPSD} imply that we must have
\begin{equation}\label{PSDtest}
\PSD(A,s),\PSD(B,s)\le 2\ell+2,\quad s=1,2,\ldots,\ell-1,
\end{equation}
which is an analog of the PSD-test for binary Legendre pairs introduced in~\cite{fgs01}. Criterion~\eqref{PSDtest} allows us to discard many sequences because they cannot be part of a Legendre pair.
(Note that we have $\PSD(A,\ell-s)=\PSD(\overline{A},s)$ which for binary sequences reduces the PSD-test $(\ref{PSDtest})$ to  $s=1,2,\ldots,\frac{\ell-1}2$. However, for quaternary sequences we need to check all $s=1,2,\ldots,\ell-1$.) 

\subsection{Norms in cyclotomic fields}

Let 
$$\mathbb{Q}(\xi_n)=\bigl\{a(\xi_n)=a_0+a_1\xi_n+\dots+a_{\varphi(n)-1}\xi_n^{\varphi(n)-1}: a_0,a_1,\dots,a_{\varphi(n)-1}\in \mathbb{Q}\bigr\}$$
denote the $n$th {\em cyclotomic field}, 
where $\varphi(n)$ is Euler's totient function.
For background on Galois theory and cyclotomic fields see for example~\cite{l06,w97}.

\subsubsection*{A condition for integers of the form $|\alpha|^2$ for some $\alpha\in \mathbb{Q}(\xi_n)$}

We recall the following well-known condition, see for example \cite[Theorem~3.1]{b88} which follows from \cite[Lemma~3]{t65}.
\begin{lemma}
\label{brock}
Suppose that $\alpha\in \mathbb{Q}(\xi_n)$ satisfies $|\alpha|^2=m$ for an integer $m$, and that $p$ is a rational prime dividing the square-free part of $m$. Denote by $f$ the order of $p$ modulo $n$ if $\gcd(p,n)=1$. Then either   
\begin{itemize}
\item $f$ is odd or
\item $f$ is even and $p^{f/2}\not\equiv -1\bmod n$.
\end{itemize}
\end{lemma}

\noindent
For $n=4$ we get the sum of two squares theorem.
\begin{lemma}\label{twosquare}
Let $m=|\alpha|^2\in\mathbb{Z}$ for some $\alpha\in \mathbb{Z}[i]$, then every rational prime  divisor~$p$ of the square-free part of $m$ is either $p=2$ or $p\equiv 1\bmod 4$.
\end{lemma}

\noindent
For $n=q>2$ a prime, Lemma~\ref{brock} simplifies to the following condition.
\begin{lemma}\label{brockp} 
Let $q>2$ be a prime.
If $m=|\alpha|^2$ is an integer for some $\alpha\in \mathbb{Q}(\xi_q)$, then for every rational prime divisor $p$ of the square-free part of $m$ we have either $p=q$ or the order of $p$ mod $q$ is odd.
\end{lemma}
For example, if $q\equiv 3\bmod 4$, then the eligible $p$ with $p\not=q$ are equivalent to a quadratic residue modulo $q$, and if $q=2^s+1$ is a Fermat prime, then either $p=q$ or $p\equiv 1\bmod q$.

\subsubsection*{Norms and automorphisms in cyclotomic fields}

The {\em automorphisms} $\psi_j$ of $\mathbb{Q}(\xi_n)$ are 
$$\psi_j(a(\xi_n))=a(\xi_n^j),\quad \gcd(j,n)=1,$$
and the group~$G_n$ of automorphisms is isomorphic to the multiplicative group~$\mathbb{Z}_n^*$ of the residue class ring $\mathbb{Z}_n$ modulo~$n$.

The {\em absolute norm} of $a(\xi_n)$ is
$$N(a(\xi_n))=\prod_{\psi \in G_n} \psi(a(\xi_n))=\prod_{\textstyle\genfrac{}{}{0pt}{}{j=1}{\gcd(j,n)=1}}^{n-1}\kern-10pt a(\xi_n^j)\ \in \mathbb{Q}$$
and if $a(\xi_n)\in \mathbb{Z}[\xi_n]$ we get $N(a(\xi_n))\in \mathbb{Z}.$

By the {\em Fundamental Theorem of Galois Theory} there is a bijection between the subgroups $U$ of $G_n$ and the subfields $\mathbb{F}$ of $\mathbb{Q}(\xi_n)$:
$$U\mapsto \mathbb{F}=\Fix(U)=\{\xi\in \mathbb{Q}(\xi_n): \psi(\xi)=\xi \text{ for all } \psi \in U\}$$
and $$\mathbb{F}\mapsto U=\Aut(\mathbb{F})=\{\psi\in G_n : \psi(\xi)=\xi \text{ for all }\xi\in \mathbb{F}\},$$
respectively.
The {\em relative norm} of $a(\xi_n)\in\Q(\xi_n)$ in $\Fix(U)$ is 
$$N_U(a(\xi_n))=\prod_{\psi \in U} \psi(a(\xi_n))\ \in \Fix(U).$$
In particular, if $n\equiv 0\bmod 4$, then $\mathbb{Q}(i)$ is a subfield of $\mathbb{Q}(\xi_n)$
and for 
$$U_4=\{\psi_{4j+1}: \gcd(4j+1,n)=1\}$$ we get
\begin{equation}\label{NU4}
N_{U_4}(a(\xi_n))=\prod_{\textstyle\genfrac{}{}{0pt}{}{j=0}{\gcd(4j+1,n)=1}}^{n/4-1} \kern-16pt a(\xi_n^{4j+1})\ \in \mathbb{Q}(i),
\end{equation}
since $i^{4j+1}=i$.
Moreover, if some prime $p>2$ divides~$n$, then $\mathbb{Q}(\xi_p)$ is a subfield of 
$\mathbb{Q}(\xi_n)$ and for 
$$U_p=\{\psi_{jp+1} : \gcd(jp+1,n)=1\}$$
we get
$$N_{U_p}(a(\xi_n))=\prod_{\textstyle\genfrac{}{}{0pt}{}{j=0}{\gcd(jp+1,n)=1}}^{n/p-1}\kern-16pt a(\xi_n^{jp+1})\ \in \mathbb{Q}(\xi_p),$$
since $\xi_p^{jp+1}=\xi_p$.

We have $\psi_{n-1}\not\in U_4$ and $\psi_{n-1}\not\in U_p$ for any prime divisor $p>2$ of $n$ since 
$n-1\equiv -1\bmod 4p$. Note that 
$$\psi_{n-1}(a(\xi_n))=\overline{a(\xi_n)},\quad \mbox{and}\quad a(\xi_n)\psi_{n-1}(a(\xi_n))=|a(\xi_n)|^2,$$ and by the transitivity of the norm we get
$$N(a(\xi_n))=\bigl|N_{U_4}(a(\xi_n))\bigr|^2.$$
In addition, write 
$$b(\xi_p)=N_{U_p}(a(\xi_n))\in \Q(\xi_p)$$
to obtain
$$N(a(\xi_n))=\left|\prod_{j=1}^{(p-1)/2}b(\xi_p^j) \right|^2,$$
that is, the absolute norm of $b(\xi_p)\in \Q(\xi_p)$ equals the absolute norm of $a(\xi_n)\in\Q(\xi_n)$, again by the transitivity of the norm.
Hence, for any $a(\xi_n)\in \mathbb{Z}[\xi_n]$ and for any prime divisor $p>2$ of~$n$, we have
\begin{equation}\label{alphap}
N(a(\xi_n))=|\alpha_4|^2=|\alpha_p|^2\in \mathbb{Z}
\end{equation}
for some $\alpha_4\in \mathbb{Q}(i)$ and some $\alpha_p\in \mathbb{Q}(\xi_p)$.

\section{General ideas for even \texorpdfstring{$\ell$}{l}}
\label{general}
Let $\ell$ be even. 
Conditions on $\PSD\left(A,\frac{\ell}{2}\right)$ and $\PSD\left(B,\frac{\ell}{2}\right)$ can be used to separate the search for $[a_0,a_2,\ldots,a_{\ell-2}]$ and $[a_1,a_3,\ldots,a_{\ell-1}]$, or $[b_0,b_2,\ldots,b_{\ell-2}]$ and $[b_1,b_3,\ldots,b_{\ell-1}]$, respectively. More precisely, put $k=\frac{\ell}{2}$, and let
\begin{equation}\label{evenodd1}
\alpha_0=\sum_{j=0}^{k-1} a_{2j},\qquad \alpha_1=\sum_{j=0}^{k-1}a_{2j+1},
\end{equation}
\begin{equation}\label{evenodd2}
\beta_0=\sum_{j=0}^{k-1} b_{2j},\qquad \beta_1=\sum_{j=0}^{k-1}b_{2j+1}.
\end{equation}
Obviously, by~\eqref{balance} we have
$$\alpha=\alpha_0+\alpha_1=0,\qquad \beta=\beta_0+\beta_1=1+i.$$
We get
$$
\PSD(A,k)=\left|\sum_{j=0}^{\ell-1}a_j (-1)^j\right|^2=|\alpha_0-\alpha_1|^2=4|\alpha_0|^2=4|\alpha_1|^2$$
and
$$
\PSD(B,k)=|\beta_0-\beta_1|^2=|2\beta_0-1-i|^2=|2\beta_1-1-i|^2.$$
In particular, 
\begin{equation}\label{A_2}
 \PSD(A,k) \quad\mbox{and}\quad  \PSD(B,k)\quad \mbox{are sums of two integer squares}. 
\end{equation}

Let $N_b$ be the number of elements $b\in \{-1,+1,-i,+i\}$ in the subsequence $[a_0,a_2,\ldots,a_{\ell-2}]$. 
Then
\begin{equation}\label{alpha0iskmod2}
  |\alpha_0|^2=(N_1-N_{-1})^2+(N_i-N_{-i})^2\equiv N_1+N_{-1}+N_i+N_{-i}\equiv k \bmod 2
\end{equation}
and thus 
\begin{equation}\label{4k}\PSD(A,k)\equiv 4k\equiv 2\ell \bmod 8.
\end{equation}
The conditions $(\ref{sumPSD})$, $(\ref{A_2})$, and $(\ref{4k})$ can now be used to construct a list of eligible pairs $\bigl(\PSD(A,k),\PSD(B,k)\bigr)$. This allows us to
separate the search for $[a_0,a_2,\ldots,a_{\ell-2}]$ and $[a_1,a_3,\ldots,a_{\ell-1}]$ as well as 
for $[b_0,b_2,\ldots,b_{\ell-2}]$ and $[b_1,b_3,\ldots,b_{\ell-1}]$.

\begin{example}\label{ex4}
For $\ell=6$ there are two pairs
\[
  \bigl(\PSD(A,3),\PSD(B,3)\bigr)=\bigl(\PSD(A,3),14-\PSD(A,3)\bigr)
\]
with $\PSD(A,3)\equiv 4\bmod 8$, namely $(4,10)$ and $(12,2)$. However, since $12$ is not the sum of two integer squares we must have 
$$\bigl(\PSD(A,3),\PSD(B,3)\bigr)=(4,10).$$
\end{example}

From now on, we assume that $k$ is even, i.e., that $\ell$ is divisible by~$4$. In this case 
\[
  \DFT\biggl(A,\frac{\ell}{4}\biggr),\ 
  \DFT\biggl(A,\frac{3\ell}{4}\biggr),\
  \DFT\biggl(B,\frac{\ell}{4}\biggr),\
  \DFT\biggl(B,\frac{3\ell}{4}\biggr)\in \mathbb{Z}[i],
\]
which means that
\begin{equation}\label{PSD14}
  \PSD\biggl(A,\frac{\ell}{4}\biggr),\
  \PSD\biggl(A,\frac{3\ell}{4}\biggr),\ 
  \PSD\biggl(B,\frac{\ell}{4}\biggr),\ 
  \PSD\biggl(B,\frac{3\ell}{4}\biggr)
\end{equation}
 are all sums of two integer squares.

Now put 
$$\alpha_j'=\sum_{s=0}^{\ell/4-1}a_{4s+j},\quad j=0,1,2,3.$$
Note that 
$$\alpha_0'+\alpha_1'+\alpha_2'+\alpha_3'=\alpha=0,$$
and
\begin{align*}
  \alpha_0'+\alpha_2' &= \alpha_0, \\
  \alpha_1'+\alpha_3' &= \alpha_1=-\alpha_0.
\end{align*}
We get
\begin{align*}
    \PSD\biggl(A,\frac{\ell}{4}\biggr) &= \bigl|\alpha_0'-\alpha_2'+(\alpha_1'-\alpha_3')i\bigr|^2=\bigl|2\alpha_0'-\alpha_0+(2\alpha_1'+\alpha_0)i\bigr|^2 \\
    &\equiv 2|\alpha_0|^2\equiv 0 \bmod 4
\end{align*}
by employing~\eqref{alpha0iskmod2}, and by using the identity
$$|x+2y|^2=|x|^2+4|y|^2+4 \Re(x\overline{y})\quad \mbox{with}\quad x=(i-1)\alpha_0.$$
Since we have
$$\PSD\biggl(B,\frac{\ell}{4}\biggr)=2\ell+2-\PSD\biggl(A,\frac{\ell}{4}\biggr)\equiv 2-\PSD\biggl(A,\frac{\ell}{4}\biggr)\bmod 8,$$
which in particular implies that $\PSD\bigl(B,\frac{\ell}{4}\bigr)$ is even,
and since $\frac12\PSD\bigl(B,\frac{\ell}{4}\bigr)$ is a sum of two squares, which by Lemma~\ref{twosquare} is impossible if $\PSD\bigl(B,\frac{\ell}{4}\bigr)\equiv 6 \bmod 8$, we get
\begin{equation}\label{0mod8}
\PSD\biggl(A,\frac{\ell}{4}\biggr)\equiv 0\bmod 8.
\end{equation}

Write $\alpha_0=u+iv$ with integers $u$ and $v$. We will show that $u\equiv v\equiv 1\bmod 2$ is not possible.
Assume to the contrary that $u$ and $v$ are both odd. Write 
$$\alpha_0'=x_0+iy_0, \quad \alpha_1'=x_1+iy_1,\quad x_0,y_0,x_1,y_1\in \mathbb{Z},$$
and observe
$$x_0+y_0+x_1+y_1\equiv \frac{\ell}{2}\equiv 0\bmod 2.$$
Hence, 
\begin{align*}
    \PSD\biggl(A,\frac{\ell}{4}\biggr) &= (2x_0-u-2y_1-v)^2+(2y_0-v+2x_1+u)^2\\
    &= 4\bigl((x_0-y_1)^2+(y_0+x_1)^2\bigr)
    -4\bigr((x_0-y_1)(u+v)-(x_1+y_0)(u-v)\bigl)\\
    &\quad +(u+v)^2+(u-v)^2.
\end{align*}
Since
$(x_0-y_1)^2+(y_0+x_1)^2\equiv x_0+y_1+y_0+x_1\equiv 0\bmod 2$, and $u\pm v$ is even, we get
$$\PSD\biggl(A,\frac{\ell}{4}\biggr)\equiv(u+v)^2+(u-v)^2\equiv 2(u^2+v^2)\bmod 8.$$
But our assumption that $u$ and $v$ are odd, that is $u^2\equiv v^2\equiv 1\bmod 8$, implies that $\PSD\bigl(A,\frac{\ell}{4}\bigr)\equiv 4\bmod 8$
which contradicts~\eqref{0mod8}. Therefore, and because $u+v\equiv \frac{\ell}{2}\equiv 0\bmod 2$, we obtain $u\equiv v\equiv 0\bmod 2$, and consequently
\begin{equation}\label{even}
\PSD\biggl(A,\frac{\ell}{2}\biggr)=4|\alpha_0|^2=4(u^2+v^2)\equiv 0\bmod 16.
\end{equation}
(Note that $(\ref{0mod8})$ is equivalent to \cite[Theorem 3.2]{kw24}. However, the proof in \cite{kw24} was incomplete which is fixed now.) \\

\begin{example}
For $\ell=32$ the pairs $\bigl(\PSD\bigl(A,\frac{\ell}{2}\bigr),\PSD\bigl(B,\frac{\ell}{2}\bigr)\bigr)$
with $(\ref{even})$ are $(0,66),(16,50),(32,34),(48,18),(64,2)$. However, since $48$ and $66$ are not sums of two squares, $(0,66)$ and $(48,18)$ are not eligible.\\
The eligible pairs $\bigl(\PSD\bigl(A,8\bigr),\PSD\bigl(B,8\bigr)\bigr)$ and $\bigl(\PSD\bigl(A,24\bigr),\PSD\bigl(B,24\bigr)\bigr)$ are given in Section~\ref{constructions}.
\end{example}

We state the results of this section as a proposition.
\begin{proposition}
Let $\ell$ be even and let $(A,B)$ be a quaternary Legendre pair of length~$\ell$. Then the following three assertions hold:
\begin{enumerate}
\item We have the following condition on $\left(\PSD\bigl(A,\frac{\ell}{2}\bigr),\PSD\bigl(B,\frac{\ell}{2}\bigr)\right)$:  
$$\PSD\bigl(A,\tfrac{\ell}{2}\bigr)\equiv 2\ell \bmod 8, \quad \PSD\bigl(B,\tfrac{\ell}{2}\bigr)=2\ell+2-\PSD\bigl(A,\tfrac{\ell}{2}\bigr),$$
and both square-free parts of $\PSD\bigl(A,\frac{\ell}{2}\bigr)$ and $\PSD\bigl(B,\frac{\ell}{2}\bigr)$ are not divisible by a rational prime $p\equiv 3\bmod 4$. 

\item The pairs $(\alpha_0,\alpha_1)$ and $(\beta_0,\beta_1)$, defined by \eqref{evenodd1} and \eqref{evenodd2}, satisfy
$$\PSD\bigl(A,\tfrac{\ell}{2}\bigr)=4|\alpha_0|^2=4|\alpha_1|^2$$ 
and
$$\PSD\bigl(B,\tfrac{\ell}{2}\bigr)=|2\beta_0-1-i|^2=|2\beta_1-1-i|^2.$$

\item If $\ell\equiv 0\bmod 4$, then, in addition, we have\\
$$\PSD\bigl(A,\tfrac{\ell}{2}\bigr)\equiv 0 \bmod 16,$$
$$\PSD\bigl(A,\tfrac{\ell}{4}\bigr)\equiv 0\bmod 8,$$
$$\PSD\bigl(B,\tfrac{\ell}{4}\bigr)=2\ell+2-\PSD\bigl(A,\tfrac{\ell}{4}\bigr)$$
and both square-free parts of $\PSD\bigl(A,\frac{\ell}{4}\bigr)$ and $\PSD\bigl(B,\frac{\ell}{4}\bigr)$ are not divisible by a rational prime $p\equiv 3\bmod 4$.
\end{enumerate}
\end{proposition}

\section{Condition for certain products of PSDs}
\label{prod}

Let $d>1$ be a divisor of $\ell$ and let $n=\lcm(4,d)$.
Then
$$\DFT\biggl(A,\frac{\ell}{d}\biggr)\in \mathbb{Z}[\xi_n].$$
%
%
%
By $(\ref{NU4})$, $(\ref{alphap})$ and 
$$\PSD\left(A,\frac{\ell}{d}\right)=\left|\DFT\left(A,\frac{\ell}{d}\right)\right|^2$$
we have the absolute norm from $\mathbb{Q}(\xi_n)$ of the DFTs
$$N\biggl(\DFT\biggl(A,s\frac{\ell}{d}\biggr)\biggr)=\prod_{\textstyle\genfrac{}{}{0pt}{}{j=0}{\gcd(4j+s,n)=1}}^{n/4-1}\kern-14pt\PSD\biggl(A,(4j+s\bmod d)\frac{\ell}{d}\biggr),\quad s\in \{1,3\}.$$
Note that for $d\not\equiv 0\bmod 4$ we get
$$N\biggl(\DFT\biggl(A,\frac{\ell}{d}\biggr)\biggr)=N\biggl(\DFT\biggl(A,3\frac{\ell}{d}\biggr)\biggr)$$
since $4j+1\equiv 3\bmod d$ for some $j$. For $d\equiv 0\bmod 4$ there is no solution $j$ of $4j+1\equiv 3\bmod d$ and the factors of these products use different arguments of the DFTs.

By $(\ref{alphap})$, Lemma~\ref{twosquare} and Lemma~\ref{brockp} are applicable and we get the following result for products of PSDs, that is, absolute norms of DFTs.

\begin{theorem}\label{prodcond}
Let $d>1$ be any divisor of $\ell$ and let $n=\lcm(4,d)$.
Then any prime divisor $p$ of the square-free part of  
\begin{eqnarray*}    
\prod_{\psi\in G_n}\psi\left(\DFT\left(A,\frac{\ell}{d}\right)\right)&=&\prod_{\psi \in U_4}\psi\left(\PSD\left(A,\frac{\ell}{d}\right)\right)\\&=&\prod_{\textstyle\genfrac{}{}{0pt}{}{j=0}{\gcd(4j+1,n)=1}}^{n/4-1}\kern-14pt\PSD\biggl(A,(4j+1\bmod d)\frac{\ell}{d}\biggr)\in \mathbb{Z}
\end{eqnarray*}
is 
$$p=2 \quad \mbox{or}\quad p\equiv 1 \bmod 4$$
and for any prime divisor $q>2$ of $n$ we have
$$p=q\quad \mbox{or}\quad \mbox{the order of $p$ modulo $q$ is odd}.$$
If $d\equiv 0\bmod 4$, the same conditions are true for
\begin{eqnarray*}\prod_{\psi\in G_n}\psi\left(\DFT\left(A,3\frac{\ell}{d}\right)\right)&=& 
\prod_{\psi\in U_4}\psi\left(\PSD\left(A,3\frac{\ell}{d}\right)\right)\\
&=&\prod_{\textstyle\genfrac{}{}{0pt}{}{j=0}{\gcd(4j+3,n)=1}}^{n/4-1}\kern-14pt\PSD\biggl(A,(4j+3\bmod d)\frac{\ell}{d}\biggr)\in \mathbb{Z}.
\end{eqnarray*}
\end{theorem}

Using the Chinese remainder theorem the different conditions modulo $4$ and modulo $p$ can be combined to one condition modulo
$$4\prod\limits_{q\mid n}q,$$ 
where the product is taken over all odd prime divisors $q$ of $n$.


For any sequence $A$ the result is true but there may be no complementary~$B$ satisfying the condition for
$$\prod_{\textstyle\genfrac{}{}{0pt}{}{j=0}{\gcd(4j+s,n)}}^{n/4-1}\kern-8pt\left(2\ell+2-\PSD\biggl(A,(4j+s\bmod d)\frac{\ell}{d}\right)\biggr).$$
Note that for $d=2$ and $d=4$ we get the tests of the previous section.\\

\begin{example}
For $\ell=6$ and $d=3$ or $d=6$ we get $n=12$ and any prime divisor~$p$ of the square-free parts of
$$\PSD(A,2)\PSD(A,4)\quad \mbox{and}\quad \PSD(A,1)\PSD(A,5)$$
must satisfy both conditions
$$\mbox{either} \quad p=2\quad\mbox{or}\quad p\equiv 1\bmod 4$$
and 
$$\mbox{either} \quad p=3\quad\mbox{or}\quad p\equiv 1\bmod 3.$$
Since $2\not\equiv 1\bmod 3$ and $3\not\equiv 1\bmod 4$ and by the Chinese remainder theorem this is equivalent to 
$p\equiv 1\bmod 12$.

We have $4^6=4096$ quaternary sequences of length $6$. There are 
$100$ sequences $B$ with 
$$b_0=1\quad\mbox{and}\quad\beta=1+i.$$ By Example \ref{ex4}, we must have 
$$\PSD(B,3)=10$$ and only 
$36$ sequences $B$ remain. 
$20$ sequences $B$ pass the PSD-test \eqref{PSDtest}
and only $4$ of them also the last PSD-tests of this example.
All these $B$ have the same PAFs, $\PAF(B,1)=\PAF(B,2)=0$, $\PAF(B,3)=-4$ and there are $12$ sequences $A$ with $a_0=1$ and the desirable PAFs,
$\PAF(A,1)=\PAF(A,2)=-2$, $\PAF(A,3)=2$. 
\end{example}


%

\section{Legendre pairs for \texorpdfstring{$\ell=28,30,32$ and $34$}{l=28, 30, 32 and 34}}
\label{constructions}

In this section, we report some computational findings, namely quaternary Legendre pairs of the previously open lengths $28,30$, and $32$, that were made possible by the theoretical results presented in the previous sections. Our strategy is as follows: we fix an eligible pair $\bigl(\PSD\bigl(A,\frac{\ell}{2}\bigr),\PSD\bigl(B,\frac{\ell}{2}\bigr)\bigr)$, from which we can deduce the possible choices for $\alpha_0,\alpha_1,\beta_0,\beta_1$ (see Section~\ref{general}). Fix one such choice and generate the set~$S_0$ (resp.~$S_1$) of all $\{+1,-1,+i,-i\}$-sequences of length $\ell/2$ whose sum equals $\alpha_0$ (resp.~$\alpha_1$). Then $S_0\times S_1$ (with suitable interlacing applied to each member) forms a set of candidate $A$-sequences. The same is done analogously and independently for the $B$-sequences. Certain optimizations can be implemented, such as modding out cyclic shifts and constant multiples of sequences, and precomputing the DFTs of the subsequences in $S_0$ and $S_1$ so that the DFT of a member of $S_0\times S_1$ is obtained by a single addition.

The main point of our investigations, and the reason why such searches become feasible at all, is the fact that many of these candidate sequences can be filtered out at an early stage, i.e., before it is tried to find a corresponding partner sequence for a Legendre pair. In particular, we have the following tests at our disposal:
\begin{itemize}
    \item[(T1)] the PSD test \eqref{PSDtest},
    \item[(T2)] the $\PSD\bigl(\cdot,\frac{\ell}{4}\bigr)$ and $\PSD\bigl(\cdot,\frac{3\ell}{4}\bigr)$ tests~\eqref{0mod8} in the case that $\ell$ is divisible by~$4$,
    \item[(T3)] a set of product tests that arise from the results of Section~\ref{prod}.
\end{itemize}
We store only those candidate $A$- and $B$-sequences that pass all tests, and then try to find matches between them that form Legendre pairs.

\subsection{The case \texorpdfstring{$\ell=28$}{l=28}}

By  $(\ref{sumPSD})$, $(\ref{A_2})$ and $(\ref{even})$ there are only two eligible pairs $\bigl(\PSD(A,14),\PSD(B,14)\bigr)$:
$$(0,58),\ (32,26).$$
By $(\ref{sumPSD})$, $(\ref{PSD14})$, $(\ref{0mod8})$ and Lemma~\ref{twosquare} the eligible pairs $\bigl(\PSD(A,7),\PSD(B,7)\bigr)$, $\bigl(\PSD(A,21),\PSD(B,21)\bigr)$ are 
$$(0,58),\ (8,50),\ (32,26),\ (40,18).$$
We also need that the prime divisors $p$ of the square-free parts of the integers
\newcommand{\alprod}[1]{\prod_{\kern-30pt #1\kern-30pt}}
\begin{alignat*}{2}
  \kern30pt
  &  \alprod{j=1}^6\ \bigl(58-\PSD(A,4j)\bigr),
  && \alprod{j=1}^6\ \bigl(58-\PSD(B,4j)\bigr),
  \\[1ex]
  &  \alprod{\textstyle\genfrac{}{}{0pt}{}{j=0}{j\ne 3}}^6\ \bigl(58-\PSD(A,4j+2)\bigr),
  \kern30pt
  && \alprod{\textstyle\genfrac{}{}{0pt}{}{j=0}{j\ne 3}}^6\ \bigl(58-\PSD(B,4j+2)\bigr),
  \\[1ex]
  &  \alprod{j\in \{1,5,9,13,17,25\}}\ \bigl(58-\PSD(A,j)\bigr),
  && \alprod{j\in \{1,5,9,13,17,25\}}\ \bigl(58-\PSD(B,j)\bigr),
  \\[1ex]
  &  \alprod{j\in \{3,11,15,19,23,27\}}\ \bigl(58-\PSD(A,j)\bigr),
  && \alprod{j\in \{3,11,15,19,23,27\}}\ \bigl(58-\PSD(B,j)\bigr)
\end{alignat*}
are either $p=2$ or $p\equiv 1,9$ or $25\bmod 28$.

We performed an exhaustive search for Legendre pairs of length~$28$, which means for $\bigl(\PSD(A,14),\PSD(B,14)\bigr)=(0,58)$ to generate 9,909,733,287,168 (resp.~5,164,090,709,778) candidates for the $A$- (resp.~$B$-) sequences, while for $\bigl(\PSD(A,14),\PSD(B,14)\bigr)=(32,26)$ it means to generate 3,372,055,150,464 (resp.~11,618,309,811,108) candidates for the $A$- (resp.~$B$-) sequences. To demonstrate the effectiveness of our tests, we look at the $A$-sequences of the latter case: after applying tests (T1) and (T2) only 1,028,107,232 sequences survive ($\sim 0.03\%$). This number can be further reduced by applying test~(T3), after which only 289,305,112 sequences remain, that is, a further reduction by a factor 3.5 approximately. In total, we found 529,152 qLPs with $\PSD(A,14)=0$ and 383,328 qLPs with $\PSD(A,14)=32$. All of them satisfy the predicted properties from Theorem~\ref{prodcond}, of course. The total time to perform this search was about 46 CPU days.

\begin{example}
We get the Legendre pair 
\begin{align*}
A &= [-1, -1, -i, -1, -i, -i, -i, 1, i, -i, -1, 1, i, 1,\\
&\qquad {-i}, -1, -1, i, i, i, -i, i, 1, 1, i, i, -i, 1],\\
B &= [-1, i, -i, 1, i, -i, i, i, i, -1, -i, 1, -i, -1,\\
&\qquad i, i, -i, i, -i, 1, -i, 1, 1, -1, -i, -1, i, 1],
\end{align*}
with the following PSD values at multiples of $\frac{\ell}{4}$:
\begin{align*}
(\PSD(A,14),\PSD(B,14)) &= (32,26), \\
(\PSD(A,7),\PSD(B,7)) &= (8,50), \\
(\PSD(A,21),\PSD(B,21)) &= (40,18).
\end{align*}

Moreover, we can verify that all PSD products satisfy our conditions:
\begin{align*}
& \alprod{j=1}^6\ \PSD(A,4j) = 164204096=2^6\cdot 7^2\cdot 52361, \\
& \alprod{j=1}^6\ \PSD(B,4j) = 315341888=2^6\cdot 1933 \cdot 2549, \\
& \alprod{\textstyle\genfrac{}{}{0pt}{}{j=0}{j\not=3}}^6\ \PSD(A,4j+2) =
  340963904 = 2^6 \cdot 29 \cdot 183709, \\
& \alprod{\textstyle\genfrac{}{}{0pt}{}{j=0}{j\not=3}}^6\ \PSD(B,4j+2) =
  20892992 = 2^6 \cdot 29 \cdot 11257, \\[1ex]
& \alprod{j\in \{1,5,9,13,17,25\}}\ \PSD(A,j) =
  120847168 = 2^6 \cdot 13^2 \cdot 11173, \\[1ex]
& \alprod{j\in \{1,5,9,13,17,25\}}\ \PSD(B,j) =
  510607168 = 2^6 \cdot 7978237, \\[1ex]
& \alprod{j\in \{3,11,15,19,23,27\}}\ \PSD(A,j) =
  288564032 = 2^6 \cdot 113 \cdot 39901, \\[1ex]
& \alprod{j\in \{3,11,15,19,23,27\}}\ \PSD(B,j) =
  36675136 = 2^6 \cdot 757^2.
\end{align*}
Note that all prime numbers~$p$ appearing in the square-free parts of the above products satisfy $p\equiv1\bmod28$.

\end{example}

\subsection{The case \texorpdfstring{$\ell=30$}{l=30}}

There are only three eligible pairs for $\bigl(\PSD(A,15),\PSD(B,15)\bigr)$:
$$(4,58),\ (36,26),\ (52,10).$$
Moreover, we get the following conditions:
\begin{enumerate}
\item Any prime divisor $p$ of the square-free parts of 
\begin{alignat*}{2}
  & \PSD(A,10)\PSD(A,20), &\quad& \PSD(B,10)\PSD(B,20), \\
  & \PSD(A,5)\PSD(A,25), && \PSD(B,5)\PSD(B,25)
\end{alignat*}
satisfies $p\equiv 1\bmod 12$.

\item Any prime divisor $p$ of the square-free parts of
\begin{align*}
  & \PSD(A,6)\PSD(A,12)\PSD(A,18)\PSD(A,24), \\
  & \PSD(B,6)\PSD(B,12)\PSD(B,18)\PSD(B,24), \\
  & \PSD(A,3)\PSD(A,9)\PSD(A,21)\PSD(A,27), \\
  & \PSD(B,3)\PSD(B,9)\PSD(B,21)\PSD(B,27)
\end{align*}
is $p=5$ or satisfies $p\equiv 1\bmod 20$.

\item Any prime divisor $p$ of the square-free parts of
\begin{alignat*}{2}
 &  \alprod{\textstyle\genfrac{}{}{0pt}{}{j=1}{\gcd(j,15)=1}}^{14}\ \PSD(A,2j),
 && \alprod{\textstyle\genfrac{}{}{0pt}{}{j=1}{\gcd(j,15)=1}}^{14}\ \PSD(B,2j),
 \\[1ex]
 &  \alprod{\textstyle\genfrac{}{}{0pt}{}{j=0}{\gcd(2j+1,15)=1}}^{14}\ \PSD(A,2j+1),
 \kern30pt
 && \alprod{\textstyle\genfrac{}{}{0pt}{}{j=0}{\gcd(2j+1,15)=1}}^{14}\ \PSD(B,2j+1)
\end{alignat*}
satisfies $p\equiv 1\bmod 60$.
\end{enumerate}


Because of the computational complexity, we did not perform an exhaustive search for $\ell=30$, but we stopped the computations as soon as some qLPs were found (that is, after about 65 hours in each of the three cases). We obtained 71 qLPs with $\PSD(A,15)=4$, then 111 qLPs with $\PSD(A,15)=36$, and finally 93 qLPs with $\PSD(A,15)=52$. Here are two quaternary Legendre pairs for $\ell=30$: 
\begin{align*}
A &= [-1, i, -i, 1, -i, i, -i, -i, -i, i, i, -i, 1, -1, -i,\\
&\qquad i, i, 1, -i, -i, -i, -i, i, i, i, -1, i, i, i, -i],\\
B &= [-1, 1, -i, -i, i, -i, -i, i, -i, -1, i, 1, i, 1, i,\\
&\qquad i, -1, -i, i, -i, i, 1, i, 1, i, -1, -i, -i, 1, -1],
\end{align*}
and
\begin{align*}
A &= [-i, -i, -1, i, -1, 1, i, -1, i, 1, 1, i, i, i, -1,\\
&\qquad i, -i, -1, -i, 1, 1, -i, -1, -i, -1, 1, i, -i, -i, 1],\\
B &= [-1, -i, i, 1, -1, -1, i, -i, 1, 1, 1, -i, i, -1, -1,\\
&\qquad 1, i, -i, -1, -1, 1, i, -i, 1, i, 1, -i, i, 1, -1].
\end{align*}

A different Legendre pair of length $30$ is given by Jedwab's and Pender's general construction \cite{jp24}:
\begin{align*}
A &= \left[\left(\frac{2-1}{61}\right),\left(\frac{2\cdot 4-1}{61}\right),\ldots,\left(\frac{2\cdot 4^{29}-1}{61}\right)\right],\\
B &= \left[i,\left(\frac{4-1}{61}\right),\left(\frac{4^2-1}{61}\right),\ldots,\left(\frac{4^{29}-1}{61}\right)\right],
\end{align*}
where $\left(\frac{.}{.}\right)$ is the Legendre symbol.

\subsection{The case \texorpdfstring{$\ell=32$}{l=32}}

The eligible values for
$\bigl(\PSD(A,16),\PSD(B,16)\bigr)$ are 
$$(16,50),\ (32,34),\ (64,2),$$
while for
$\bigl(\PSD(A,8),(\PSD(B,8)\bigr)$ and $\bigl(\PSD(A,24),\PSD(B,24)\bigr)$ the eligible values are 
$$(8,58),\ (16,50),\ (32,34),\ (40,26),\ (64,2).$$
Moreover, the integers
\begin{alignat*}{2}
  &  \PSD(A,4)\PSD(A,20),
  && \PSD(B,4)\PSD(B,20), \\
  &  \PSD(A,12)\PSD(A,28), \quad
  && \PSD(B,12),\PSD(B,28), \\
  &  \alprod{j=0}^3\ \PSD(A,8j+2),
  && \alprod{j=0}^3\ \PSD(B,8j+2), \\
  &  \alprod{j=0}^3\ \PSD(A,8j+6),
  && \alprod{j=0}^3\ \PSD(B,8j+6), \\
  &  \alprod{j=0}^7\ \PSD(A,4j+1),
  && \alprod{j=0}^7\ \PSD(B,4j+1), \\
  &  \alprod{j=0}^7\ \PSD(A,4j+3),
  && \alprod{j=0}^7\ \PSD(B,4j+3)
\end{alignat*}
are all sums of two squares and thus their square-free parts are not divisible by any prime $p\equiv 3\bmod 4$.

Similar to $\ell=30$ we did not perform an exhaustive search. In total, we identified 12 qLPs of length~$32$. Here are two quaternary Legendre pairs of length $\ell=32$ with $\PSD(A,16)=32$ and $\PSD(A,16)=64$, respectively:,
\begin{align*}
A &= [-1, -i, -i, -1, i, i, -1, -1, i, i, i, 1, -i, -i, i, i,\\
&\qquad {-i}, i, -i, -i, -i, -i, -i, 1, i, 1, -i, 1, i, i, -i, i],\\
B &= [-1, i, -i, 1, i, -i, i, -i, -i, i, -i, -i, i, i, i, -i,\\
&\qquad i, -i, -1, 1, i, -1, i, -i, i, 1, -i, -i, -i, i, 1, i],
\end{align*}
and
\begin{align*}
A &= [1, -i, i, i, -1, i, i, 1, 1, -i, -1, i, -1, -i, 1, 1,\\
&\qquad {-1}, i, -i, 1, -1, 1, -1, 1, i, -i, -1, i, -i, -1, -i, -i],\\
B &= [i, 1, 1, 1, i, 1, 1, -1, -i, -1, i, -1, -1, 1, i, -i,\\
&\qquad {-i}, -1, -1, 1, -i, -1, -i, i, -i, 1, i, i, -1, -1, 1, 1].
\end{align*}

\subsection{The case \texorpdfstring{$\ell=34$}{l=34}}
We have the following PSD-conditions,
$$(\PSD(A,17),\PSD(B,17))\in \{(20,50),(36,34),(52,18),(68,2)\}$$
and any prime divisor $p$ of the square-free parts of the integers
$$\prod_{j=1}^{16} (70-\PSD(A,2j)),\prod_{j=1}^{16} (70-\PSD(B,2j)),$$
$$\prod_{\textstyle\genfrac{}{}{0pt}{}{j=0}{j\ne 8}}^{16} (70-\PSD(A,2j+1)),\prod_{\textstyle\genfrac{}{}{0pt}{}{j=0}{j\ne 8}}^{16} (70-\PSD(B,2j+1)),$$
satisfies $p=17$ or $p\equiv 1\bmod 68$. 
We found $22$ quaternary Legendre pairs of length $34$, for example,
\begin{eqnarray*}A&=&[-1, i, -1, -1, -i, 1, 1, i, 1, i, -i, -i, i, -i, -i, -i, -i,\\
&& -1, 1, 1, i, -i, -i, i, -i, i, -i, i, i, i, -1, i, i, -i],\\ 
B&=&[-1, -i, -1, i, 1, i, i, -i, -1, -i, i, 1, 1, -i, -i, i, -1,\\
&&i, -i, 1, i, i, -i, -1,
    i, i, -i, i, -i, 1, 1, -i, -i, i],
\end{eqnarray*}
with
$$(\alpha_0,\beta_0)=(-3i,-1-2i)$$
and thus
$$(\PSD(A,17),\PSD(B,17))=(36,34).$$

\section*{Acknowledgments}
 Christoph Koutschan was supported by the Austrian Science Fund (FWF): grant 10.55776/I6130.\\
 Arne Winterhof was funded in part by the Austrian Science Fund (FWF): grant 10.55776/PAT4719224.\\
 The authors would like to thank Jonathan Jedwab and Thomas Pender for pointing to their parallel and independent work \cite{jp24}.\\
 The authors also would like to thank the anonymous referees for many very useful comments.\\
 We dedicate this work to the memory of Douglas Adams.

\end{document}